\documentclass[10pt]{article}
\usepackage{eurosym}
\usepackage{amsfonts}
\usepackage{amssymb}
\usepackage{amsmath}
\usepackage{mitpress}
\usepackage{graphicx}

\newdimen\dummy
\dummy=\oddsidemargin
\input{tcilatex}

\begin{document}

\title{Comments on "Necessary optimality conditions of an optimization
problem governed by a double phase PDE"}
\author{}
\maketitle

\ \ \ \ \ \ \ \ \ \ \ \ \ \ \ \ \ \ \ \ \ \ \ \ \ \ \ \ \ \ \ \ \textbf{Omar
Benslimane and Nazih Abderrazzak Gadhi }\footnote{%
\noindent {\footnotesize LAMA, FSDM, Sidi Mohamed Ben Abdellah University,
Fes, Morocco. Emails: abderrazzak.gadhinazih@usmba.ac.ma \ \ \ \ \ , \ \ \ \
obenslimane22@gmail.com}}

\bigskip

\noindent \textbf{Abstract.} This note concerns the paper by Benslimane et
Gadhi (JMAA. doi: 10.1016/j.jmaa.2023.127117) where the authors established
necessary optimality conditions for an optimization problem $\left( \mathcal{%
P}\right) $ governed by a double phase partial differential equation $\left( 
\mathcal{P}_{f}\right) .$ Having noticed inconsistencies between the
investigated problem $\left( \mathcal{P}_{f}\right) $ and the weak
formulation of the equation that is associated with it, we propose some
modifications that, from our perspective, correct the forthmentioned issue.
Some careless mistakes and typos that caught our interest are also
underlined and rectified.

\bigskip

 \newpage

\section{Introduction}

\noindent The double-phase operator has piqued the interest of several
academics during the last decade \cite{bb1, 1, bb2, bb3, wi, bb4, bb5, bb6,
bb7}. Let $\Omega \subseteq \mathbb{R}^{n}$ be an open, bounded and regular
set with a smooth boundary $\partial \Omega $. Let $C_{0}^{\infty }\left(
\Omega \right) $ be the set of infinitely differentiable functions with
compact supports on $\Omega $ \cite{harjulehto2019generalized}. Let $%
L^{p}\left( \Omega \right) $ and $\mathcal{W}^{1,p}\left( \Omega \right) $
be the Lebesgue and the Sobolev spaces defined by 
\begin{equation*}
L^{p}\left( \Omega \right) =\left\{ u:\Omega \rightarrow \mathbb{R}%
/\int_{\Omega }|u(x)|^{p}\ dx<\infty \right\} 
\end{equation*}%
and 
\begin{equation*}
\mathcal{W}^{1,p}\left( \Omega \right) =\left\{ u:u\in L^{p}\left( \Omega
\right) \ \mbox{and}\ |\nabla u|\in L^{p}\left( \Omega \right) \right\} .
\end{equation*}%
In \cite{1}, Benslimane and Gadhi investigated the following optimal control
problem%
\begin{equation*}
\left( \mathcal{P}\right) :\left\{ 
\begin{array}{rcl}
\underset{u,f}{\min }\ \mathcal{E}\left( f,u\right) \ \ \ \ \ \ \ \ \ \ \ \
\ \ \ \ \ \ \ \ \ \ \ \ \ \ \ \ \ \ \ \ \ \ \  &  &  \\ 
\mbox{subject to}:u\in \psi \left( f\right) ,\ u\in \mathcal{W}%
_{0}^{1,p}\left( \Omega \right) ,\ f\in L^{q}\left( \Omega \right)  &  & 
\end{array}%
\right. 
\end{equation*}%
where, for each $f\in L^{q}\left( \Omega \right) ,\ \psi \left( f\right) $
is the set of all non-trivial solution of the following double phase problem 
\begin{equation}
\left( \mathcal{P}_{f}\right) :\left\{ 
\begin{array}{rcl}
-\mbox{div}\left( |\nabla u|^{p-2}\nabla u+\mu \left( x\right) \ |\nabla
u|^{q-2}\nabla u\right) =f\ \ \mbox{in}\ \ \Omega  &  &  \\ 
u=0\ \ \mbox{on}\ \ \partial \Omega . &  & 
\end{array}%
\right.   \label{12}
\end{equation}%
where $\mathcal{W}_{0}^{1,p}\left( \Omega \right) $ is the closure of $%
C_{0}^{\infty }\left( \Omega \right) $ in $W^{1,p}\left( \Omega \right) ,$ $%
\mu :\overline{\Omega }\rightarrow \mathbb{R}_{+}^{\ast }$ is a Lipschitz
continuous function and $\mathcal{E}:L^{q}\left( \Omega \right) \times 
\mathcal{W}_{0}^{1,p}\left( \Omega \right) \rightarrow \mathbb{R}$ is a
weakly lower continuous function, which is continuously Fr\'{e}%
chet-differentiable and bounded from below. The spaces $L^{q}\left( \Omega
\right) $ and $\mathcal{W}_{0}^{1,p}\left( \Omega \right) $ are equipped
with the norms 
\begin{equation*}
\left\Vert u\right\Vert _{L^{q}\left( \Omega \right) }=\left[ \int_{\Omega
}|u|^{q}\ dx\right] ^{\frac{1}{q}}\text{ and }\left\Vert u\right\Vert _{%
\mathcal{W}_{0}^{1,p}\left( \Omega \right) }=\left[ \sum_{i=1}^{n}\int_{%
\Omega }\left\vert \frac{\partial u\left( x\right) }{\partial x_{i}}%
\right\vert ^{p}dx\right] ^{\frac{1}{p}}.
\end{equation*}%
An element $u\in \mathcal{W}_{0}^{1,p}\left( \Omega \right) $ was said to be
a non-trivial solution of $\left( P\right) $ if 
\begin{equation}
\sum_{i=1}^{n}\int_{\Omega }\left( \bigg \vert\frac{\partial u\left(
x\right) }{\partial x_{i}}\bigg \vert^{p-2}\frac{\partial u\left( x\right) }{%
\partial x_{i}}+\mu \left( x\right) \bigg \vert\frac{\partial u\left(
x\right) }{\partial x_{i}}\bigg \vert^{q-2}\frac{\partial u\left( x\right) }{%
\partial x_{i}}\right) \ \frac{\partial \varphi (x)}{\partial x_{i}}\
dx=\int_{\Omega }f\left( x\right) \varphi \left( x\right) \ dx  \label{33}
\end{equation}%
is satisfied for all test function $\varphi \in \mathcal{W}_{0}^{1,p}\left(
\Omega \right) .$ In order to get necessary optimality conditions for the
optimal control problem $\left( \mathcal{P}\right) ,$ the authors made use
of an energy function $\mathcal{J}:\mathcal{W}_{0}^{1,p}\left( \Omega
\right) \rightarrow \mathbb{R}$ defined by 
\begin{equation*}
\mathcal{J}\left( u\right) =\frac{1}{p}\sum_{i=1}^{n}\int_{\Omega }\bigg
\vert\frac{\partial u\left( x\right) }{\partial x_{i}}\bigg \vert^{p}\ dx+%
\frac{1}{q}\sum_{i=1}^{n}\int_{\Omega }\mu \left( x\right) \ \bigg \vert%
\frac{\partial u\left( x\right) }{\partial x_{i}}\bigg \vert^{q}\
dx-\int_{\Omega }f\left( x\right) u\left( x\right) \ dx.
\end{equation*}%
Their method involved proving first the energy function $\mathcal{J}$'s $p$%
-hyperconvexity and G\^{a}teaux differentiability. This allowed them to
conclude that the solution operator $\psi :L^{q}\left( \Omega \right)
\longrightarrow \mathcal{W}_{0}^{1,p}\left( \Omega \right) $ associated with 
$\left( \mathcal{P}_{f}\right) $ is both well-defined and G\^{a}%
teaux-differentiable. Following that, they proved that $\left( \mathcal{P}%
_{f}\right) $ admits a unique non-trivial solution $\mathbb{\psi }\left(
f\right) \in \mathcal{W}_{0}^{1,p}\left( \Omega \right) ,$ which enables
them to determine necessary optimality conditions of the optimal control
problem $\left( \mathcal{P}\right) $ of interest.

\noindent In this note, we underline the fact that energy function $\mathcal{%
J}$ is not well defined and that the weak formulation of the equation $%
\left( \ref{33}\right) $ does not correspond to the double phase problem $%
\left( \ref{12}\right) ,$ but corresponds in fact to the pseudo-double phase
problem $\left( \ref{bbg}\right) ,$ defined in the next section. Adjustments
that we believe will fix these awkward circumstances are then proposed. Some
typos and careless mistakes that drew our attention are also highlighted and
then corrected.

\section{Comments and corrections}

\noindent First and foremost, simply replacing 
\begin{equation*}
p<q\text{ and }\frac{1}{q}=\frac{1}{p}-\frac{1}{n}
\end{equation*}%
with 
\begin{equation*}
q<p\text{ and }\frac{1}{p}=\frac{1}{q}-\frac{1}{n}
\end{equation*}%
renders the energy functional $\mathcal{J}$ well-defined. Accordingly, to
retain the accuracy of the proof of \cite[Lemma 7]{1}, it is clear that 
\begin{equation*}
\mathcal{X}_{i}\left( x\right) =\frac{1}{p}\left( \left\vert \frac{\partial u%
}{\partial x_{i}}+\frac{\partial v}{\partial x_{i}}\right\vert ^{p}+\mu
\left( x\right) \left\vert \frac{\partial u}{\partial x_{i}}+\frac{\partial v%
}{\partial x_{i}}\right\vert ^{q}\right) ,\ i\in I,
\end{equation*}%
should be replaced with 
\begin{equation*}
\mathcal{X}_{i}\left( x\right) =\frac{1}{q}\left( \left\vert \frac{\partial u%
}{\partial x_{i}}+\frac{\partial v}{\partial x_{i}}\right\vert ^{p}+\mu
\left( x\right) \left\vert \frac{\partial u}{\partial x_{i}}+\frac{\partial v%
}{\partial x_{i}}\right\vert ^{q}\right) ,\ i\in I,
\end{equation*}%
within this proof. On the other hand, in order to address the unconstency
between $\left( \ref{12}\right) $ and $\left( \ref{33}\right) ,$ we propose
to replace the standard $p$-Laplacian in $\left( \mathcal{P}_{f}\right) $ by
the pseudo-$p$-Laplacian. By doing this, the optimal control problem $\left( 
\mathcal{P}\right) $ becames 
\begin{equation*}
\left( \mathcal{P}\right) :\left\{ 
\begin{array}{rcl}
\underset{u,f}{\min }\ \mathcal{E}\left( f,u\right) \ \ \ \ \ \ \ \ \ \ \ \
\ \ \ \ \ \ \ \ \ \ \ \ \ \ \ \ \ \ \ \ \ \ \  &  &  \\ 
\mbox{subject to}:u\in \psi \left( f\right) ,\ u\in \mathcal{W}%
_{0}^{1,p}\left( \Omega \right) ,\ f\in L^{q}\left( \Omega \right)  &  & 
\end{array}%
\right. 
\end{equation*}%
where, for each $f\in L^{q}\left( \Omega \right) ,\ \psi \left( f\right) $
is the set of all non-trivial solution of the following pseudo-double phase
problem 
\begin{equation}
\left( \mathcal{P}_{f}\right) :\left\{ 
\begin{array}{rcl}
-\sum_{i=1}^{n}\frac{\partial }{\partial x_{i}}\left( \left\vert \frac{%
\partial u}{\partial x_{i}}\right\vert ^{p-2}\frac{\partial u}{\partial x_{i}%
}+\mu \left( x\right) \left\vert \frac{\partial u}{\partial x_{i}}%
\right\vert ^{q-2}\frac{\partial u}{\partial x_{i}}\right) =f\ \ \mbox{in}\
\ \Omega  &  &  \\ 
u=0\ \ \mbox{on}\ \ \partial \Omega . &  & 
\end{array}%
\right.   \label{bbg}
\end{equation}%
This adjustment has the advantage of preserving the relationship, which
should tie the Euler-Lagrange equation $\left( \ref{33}\right) $ to the
pseudo-double phase problem $\left( \mathcal{P}_{f}\right) .$ It is worth
mentionning that the pseudo-$p$-Laplacian \cite{bk, lion} defined by 
\begin{equation*}
\widetilde{\Delta }_{p}u:=\sum_{i=1}^{n}\frac{\partial }{\partial x_{i}}%
\left( \left\vert \frac{\partial u}{\partial x_{i}}\right\vert ^{p-2}\frac{%
\partial u}{\partial x_{i}}\right) 
\end{equation*}%
has some different features compared to the standard $p$-Laplacian 
\begin{equation*}
\Delta _{p}u:=\func{div}\left( \left\vert \nabla u\right\vert ^{p-2}\nabla
u\right) 
\end{equation*}%
appearing in $\left( \ref{12}\right) ,$ in particular concerning regularity
results. However, it is worth mentioning that these two operators coincide
in the case $p=1,$ because 
\begin{equation*}
\widetilde{\Delta }_{p}u=\Delta _{p}u=\left( \left\vert u^{\prime
}\right\vert ^{p-2}u^{\prime }\right) ^{\prime }=\left( p-1\right)
\left\vert u^{\prime }\right\vert ^{p-2}u^{\prime \prime }
\end{equation*}%
for $C^{2}$-functions. Taking into account the previously indicated
modification, 
\begin{equation*}
-\mbox{div}\left( |\nabla u|^{p-2}\nabla u+\mu \left( x\right) \ |\nabla
u|^{q-2}\nabla u\right) 
\end{equation*}%
should be replaced with%
\begin{equation*}
-\sum_{i=1}^{n}\frac{\partial }{\partial x_{i}}\left( \bigg \vert\frac{%
\partial u}{\partial x_{i}}\bigg \vert^{p-2}\frac{\partial u}{\partial x_{i}}%
+\mu \left( x\right) \bigg \vert\frac{\partial u}{\partial x_{i}}\bigg \vert%
^{q-2}\frac{\partial u}{\partial x_{i}}\right) 
\end{equation*}%
throughout the document for all $u\in \mathcal{W}_{0}^{1,p}\left( \Omega
\right) .$

\noindent Following that, we suggest that in the definition of the $\gamma $%
-hyperconvexity \cite[Definition 1]{1}, the reel $c$ be strictly positive.
In addition, to have a factually double phase problem, the codomain of the
weight-function $\mu $ should be $\mathbb{R}_{+}$ and \cite[Assumption 1]{1}
should be revised as follows.%
\begin{equation*}
\text{There exists }\mu _{1}>0\text{ such that for all }x\in \mathcal{W}%
_{0}^{1,q}\left( \Omega \right) ,\text{ we have }0\leq \mu \left( x\right)
\leq \mu _{1}.
\end{equation*}%
It should be noted that the aforementioned changes have no effect on the
proofs or findings presented in \cite{1}, but it does allow for the energy
functional to exhibit two different kinds of growth on the $\left\{ \mu
=0\right\} $ and $\left\{ \mu >0\right\} .$

\noindent Our final comment concerns the proof of \cite[Lemma 5]{1}. After
some true steps, the authors have concluded that%
\begin{equation*}
\left( h+g\right) \left( \theta x+\left( 1-\theta \right) y\right) \leq
\theta \ \left( h+g\right) \left( x\right) +\left( 1-\theta \right) \ \left(
h+g\right) \left( y\right) -\min \left( \theta ,\left( 1-\theta \right)
\right) \ \left[ c\ \Vert x-y\Vert ^{p}+c^{\prime }\ \Vert x-y\Vert ^{q}%
\right] .
\end{equation*}%
From our perspective, the remaining of the proof should be replaced by the
following: since $c^{\prime }\ \Vert x-y\Vert ^{q}\geq 0,$ we get 
\begin{equation*}
\left( h+g\right) \left( \theta x+\left( 1-\theta \right) y\right) +c\ \min
\left( \theta ,\left( 1-\theta \right) \right) \ \Vert x-y\Vert ^{p}\leq
\theta \ \left( h+g\right) \left( x\right) +\left( 1-\theta \right) \ \left(
h+g\right) \left( y\right) ,
\end{equation*}%
which implies the $p-$hyperconvexity of the sum-function $h+g.$

\section{Conclusion}

\noindent In this note, we highlight certain inconsistencies in \cite{1} and
then propose some modifications that, from our perspective, correct them.

\end{document}